\documentclass[11pt]{article}
\usepackage{amsthm}
\usepackage{amsmath}
\usepackage{amsfonts}
\usepackage{cancel}
\usepackage{amssymb}

\usepackage[all, knot]{xy}

\newcommand{\Dorfman}[1]{\bleft  #1\bright}

\newcommand {\emptycomment}[1]{}

\newcommand{\dM}{\mathrm{d}}

\def\Z{\mathbb{Z}}
\def\R{\mathbb{R}}
\def\K{\mathbb{K}}

\def\C{\mathbb{C}}

\newcommand{\cale}{{\cal E}}
\newcommand{\calf}{{\cal F}}

\newcommand{\br}[1]{[\cdot,\cdot]}

\renewcommand{\deg}{{\rm deg}}
\newcommand{\maps}{\colon}

\newcommand{\h}{\mathfrak{h}}
\newcommand{\g}{\mathfrak{g}}

\newcommand{\calC}{{\cal C}}

\newcommand{\calL}{{\cal L}}

\newcommand{\calV}{{\cal V}}

\newcommand{\id}{\mbox{id}}

\newcommand{\ad}{\operatorname{ad}}

\newcommand{\Inn}{\operatorname{Inn}}
\newcommand{\Der}{\operatorname{Der}}

\newcommand{\End}{\operatorname{End}}

\newcommand{\Ker}{\operatorname{Ker}}

\newcommand{\glnv}{{\mathfrak g \mathfrak l}(V)}
\newcommand{\gl}{\mathfrak{gl}}
\newcommand{\half}{\textstyle{\frac{1}{2}}}
\newcommand{\three}{\textstyle{\frac{1}{3}}}
\newcommand{\four}{\textstyle{\frac{1}{4}}}
\newcommand{\six}{\textstyle{\frac{1}{6}}}
\newcommand{\eight}{\textstyle{\frac{1}{8}}}
\newcommand{\twelve}{\textstyle{\frac{1}{12}}}
\newcommand{\bleft}{[\![}
\newcommand{\bright}{]\!]}

\newcommand{\al}{\alpha}
\newcommand{\be}{\beta}
\newcommand{\ga}{\gamma}

\newcommand{\ve}{\varepsilon}

\newcommand{\om}{\omega}

\def\op{{\oplus}}

\def\la{\langle}
\def\ra{\rangle}

\newcommand{\ot}{\otimes}

\newcommand{\trr}{\triangleright}

\newcommand{\mo}{{}_{}}
\newcommand{\mi}{{}_{}}

\newcommand{\pf}{\noindent{\bf Proof.}\ }
\allowdisplaybreaks

\newtheorem{Theorem}{Theorem}[section]
\newtheorem{Proposition}[Theorem]{Proposition}
\newtheorem{Definition}[Theorem]{Definition}

\newtheorem{Example}[Theorem]{Example}
\newtheorem{Remark}[Theorem]{Remark}

\title{Lie color 2-algebras and omni-Lie color algebras}

\author{Tao Zhang}

\date{}

\begin{document}


 \maketitle

 \setcounter{section}{0}

 \vskip0.1cm

{\bf Abstract}\quad
The notions of Lie color 2-algebras  and 2-term color $L_{\infty}$-algebras over a group-graded vector space are introduced and studied.
It is proved that the category of Lie color $2$-algebras and the category of $2$-term color $L_{\infty}$-algebras are equivalent.
We construct Lie color $2$-algebras from omni-Lie color algebras and Leibniz color algebras. Some example of $\mathbb{Z}_2^2$-graded Lie color 2-algebras  are given.

\section{Introduction}

The notion of $L_{\infty}$-algebras which generalizes Lie algebras first appeared in deformation theory and then in closed string field theory, see \cite{SS,LS93}.
This structure has gained great attentions since Kontsevitch used $L_{\infty}$-morphisms to prove the existence of star-product on Poisson manifolds, see \cite{Kont}.
In \cite{Roytenberg}, it was proved that every Courant algebroid gives rise to an $L_{\infty}$-algebra.
In \cite{Rogers1,Zam},  $L_{\infty}$-algebras of Hamiltonian forms arise in multisymplectic geometry.

The algebraic theory of Lie 2-algebras was studied by Baez and Crans in \cite{Baez}.
Roughly speaking, a Lie 2-algebra is a categorification of a Lie algebra, where the underling vector space is replaced by 2-vector space and the
Jacobi identity is replaced by a nature transformation, called Jacobiator, which satisfies some coherence law.
They proved that the category of Lie  $2$-algebras and the category of $2$-term $L_{\infty}$-algebras are equivalent.
In \cite{BH}, Baez and Huerta view the four normed division algebras or more generally Poincare Lie superalgebra as Lie 2-superalgebra.
Using this idea, they explained the deep relationships between higher gauge theory and superstring theory.

Among others things, the most important examples of Lie 2-algebras are omni-Lie algebras which were introduced by Weinstein in \cite{Wei}
as a linearization of the  Courant algebroid.  An omni-Lie algebra is actually a Lie 2-algebra since
every Courant algebroid gives rise to a Lie 2-algebra.
Omni-Lie algebras are studied from several aspects and are generalized to omni-Lie algebroids and omni-Lie 2-algebras in \cite{CL, CLS, KW, SLZ}.
In a recent paper \cite{ZL}, we study dirac structures of omni-Lie superalgebras.

On the other hand,  group-graded structures appeared naturally in Algebra, Geometry and Physics.
In Algebra, quaternions, octonions and, more generally, any Clifford algebra admits $\Z_2^n$-grading \cite{AM,Lyc}.
In Geometry, the tangent bundle and cotangent bundle of a supermanifold is a $\Z_2^2$-supermanifolds, see \cite{Mol}.
In Physics,  $\Z_2^n$-grading are used in string theory and parastastical supersymmetry, see \cite{RW,Sch} and references therein.
It turns out that Courant algebroid structures on a vector bundle $E$ over a manifold  $M$  are  in one-to-one correspondence
with  integrable homological vector fields $Q$ over some graded supermanifold \cite{Roytenberg2}.
For more infirmations about Lie color algebra, see \cite{Kac,RW,Sch,SZ} and references therein.

Now a natural question arise, does there exist a categorification of a Lie color algebra?
In this paper, we give a positive answer to this question.
%How to categorify  a Lie color algebra/can we categorify Lie color algebras?
%In this paper, we introduced the concept of Lie color 2-algebras  and 2-term color $L_{\infty}$-algebras
%as a generalization of Lie 2-algebras and 2-term $L_{\infty}$-algebras.
%In other words, we study Lie color algebras in the category of 2-vector spaces.
%In the first part of this paper, we give a positive answer to this question.
The idea is very simple, if we replace the category of 2-vector spaces by the category of $G$-graded 2-vector spaces,
then a Lie color 2-algebra is actually a Lie color algebra in this category.
We also introduced the concept of 2-term color $L_{\infty}$-algebras, and proved that the category of Lie color $2$-algebras
and the category of $2$-term color $L_{\infty}$-algebras are equivalent.
%In this part, we mainly follow the work done in \cite{Baez,shengchen}.
The algebraic aspects of Lie color 2-algebras look very similar to those of Lie 2-algebras in \cite{Baez},
but this is a nontrivial generalization since they have more interesting grading.
%Their application in Geometry and Physics has still to be explained.

The second part of this paper is devoted to construct examples of  Lie color 2-algebras.
In this third subsection of Section 4, we introduce the concept of omni-Lie color algebras and construct Lie color 2-algebras from them.
This generalized the result of \cite{ZL}, but we also provided Example 4.11 which did not appeared in that paper.
Furthermore, at the last of Section 4, we give a method to construct Lie color 2-algebras from Leibniz color algebras.
Since not every Leibniz color algebra comes from omni-Lie algebra, we provide Example 4.17 to clarify this point.

The main results and organization of this paper are as follows.
In Section 2, we recall some facts about Lie color algebras.
In Section 3, we  introduce the notion of Lie color 2-algebras and 2-term color $L_\infty$-algebras.
It is proved that there is an equivalence between the category of Lie color $2$-algebras and the category of $2$-term color $L_{\infty}$-algebras.
In Section 4, we investigate some special cases of Lie color 2-algebras, such as skeletal and strict ones.
We also define omni-Lie color algebra $\mathcal{E}$ for a Lie color algebra $\g$ and a $\g$-module $V$.
At last, we construct Lie color 2-algebras from  omni-Lie color algebras and Leibniz color algebras.

\section{Preliminaries}%{Lie Color Algebras and Leibniz Color Algebras}

In this section, we will recall some facts and definitions about Lie color algebras, see \cite{Lei80,Lei11,OS,Sch,SZ} for more details.

Let $G$ be an abelian group.
A color vector space or a $G$-graded vector space is a vector space with direct sum decomposition $V=\bigoplus_{\al\in G}V_\al$.
An element of $V$ is said to be homogeneous of degree $\al$ if it is an element of $V_\al$.
We denote the degree of a homogeneous element by $|x|:=\al$ for $x\in V_{\al}$.
A sub-vector space $W$ of $V$ is graded if $W=\bigoplus_{\al\in G}(W\cap V_\al)$.

A homomorphism $ f: V \to W$ between two $G$-graded vector spaces $V$ and $W$ is a grade-preserving linear map:
$f(V_\al)\subseteq W_\al$, for all $\al\in G.$
The category of $G$-graded vector spaces is denoted by ${\rm\mathbf{Vect}}^G$.

 A map $\ve:G\times G\to \K^*$ is called a bicharacter on $G$ if the following identities hold,
 \begin{eqnarray*}
 \ve(\al, \be+\ga)&=&\ve(\al, \be)\ve(\al, \ga),\\
 \ve(\al+\be,\ga)&=&\ve(\al, \ga)\ve(\be, \ga),\\
 \ve(\al, \be)\ve(\be, \al)&=&1,
  \end{eqnarray*}
for all $\al, \be, \ga\in G$.
We assume that $G$ is a fixed abelian group and $\ve$ is a fixed bicharacter all though the paper.

\begin{Remark}
If $G=\Z_2$, $\ve(x, y)=(-1)^{xy}$, for all $x, y\in
\mathbb{Z}_2$. For more general bicharacters, see \cite{BFM}.
\end{Remark}

The concept of Lie color algebras was introduced in \cite{Sch} under the name of $\ve$-Lie algebras.

\begin{Definition}\label{def:colorLie}  A Lie color algebra or $\ve$-Lie algebras is a $G$-graded vector space $\g=
\bigoplus_{\al\in G}\g_\al $ together with a bracket
$[\cdot,\cdot]: \g \ot \g \to \g$ satisfying the following condition:

(i) graded condition: $[\g_\al ,\g_\be]\subseteq \g_{\al+\be}$, %$\forall\ \al, \be\in G$

(ii) $\ve$-skew symmetry:
\begin{align}
[x,y]+\ve(x\mi, y\mi)[y\mo, x\mo]=0,
\end{align}

(iii) $\ve$-Jacobi identity:
\begin{align}\label{J1}
J_1(x,y,z):=\ve(z\mi, x\mi)[[x\mo, y],z\mo]+\ve(x\mi, y\mi)[[y\mo,z],x\mo]+\ve(y\mi, z\mi)[[z\mo, x], y\mo]=0,
\end{align}
for all $x\in \g_{\al},y\in \g_{\be},z\in \g_{\gamma}$.
\end{Definition}
In the above definition, we write $\ve(x, y)$ instead of $\ve(|x|,|y|)=\ve(\al,\be)$.
%%%%%%%%%%%%%%%%%%%%%%%%%%%%%%%%%%%%%%%%%%%%%%%%%%%%%%%%%%%%%%%%%%%%%%%%%%%%%%%
By the $\ve$-skew symmetry condition, the $\ve$-Jacobi identity can be write in other forms:
\begin{align}\label{J2}
J_2(x,y,z)=[[x, y],z]+\ve(x\mi, y\mi + z\mi)[[y,z],x]+\ve(x\mi + y\mi, z\mi)[[z, x], y]=0,
\end{align}
and the graded Leibniz rule
\begin{align}
J_3(x,y,z)&=[x,[y,z]]-[[x,y],z]-\ve(x,y)[y,[x, z]]=0.
\end{align}
In fact, $J_3=-J_2=-\ve(x,z)J_1$.
When the bracket is not $\ve$-skew symmetric, we have the concept of Leibniz color algebras.
%%%%%%%%%%%%%%%%%%%%%%%%%%%%%%%%%%%%%%%%%%%%%%%%%%%%%%%%%%%%%%%%%%%%%%%%%%%%%%%
\begin{Definition}\label{def:colorLeibniz} A Leibniz color algebra is a $G$-graded vector space $L=\bigoplus_{\al\in G}L_
\al $ together with a bracket
$\circ: L  \times L  \to L $ satisfies $L_\al \circ L_\be\subseteq L_{\al+\be}$, $\forall \al,\be\in G$,
and the $\ve$-Leibniz rule:
\begin{align}\label{eqn:Leibniz}
x\circ (y\circ z)=(x\circ y)\circ z+\ve(x,y)\ y\circ(x\circ z),
\end{align}
for all $x\in L_{\al},y\in L_{\be},z\in L_{\gamma}$.
\end{Definition}
%%%%%%%%%%%%%%%%%%%%%%%%%%%%%%%%%%%%%%%%%%%%%%%%%%%%%%%%%%%%%%%%%%%%%%%%%%%%%%%
\begin{Remark}
For any abelian group $G$, we have the trivial bicharacter $\ve(\alpha, \beta)=1$, for all $\alpha, \beta\in G$, then this is a graded Lie algebra.
If the group $G$ is also trivial, i.e. $G=\{0\}$, then this is an ordinary Lie algebra.
\end{Remark}
%%%%%%%%%%%%%%%%%%%%%%%%%%%%%%%%%%%%%%%%%%%%%%%%%%%%%%%%%%%%%%%%%%%%%%%%%%%%%%%
\begin{Remark}
If $G=\Z_2$, $\ve(x, y)=(-1)^{xy}$, for all $x, y\in
\mathbb{Z}_2$, this is exactly Lie superalgebra \cite{Kac}.
In this case, the Jacobi identity becomes
$$(-1)^{zx}[[x, y],z]+(-1)^{xy}[[y,z],x]+(-1)^{yz}[[z, x], y]=0.$$
\end{Remark}
%%%%%%%%%%%%%%%%%%%%%%%%%%%%%%%%%%%%%%%%%%%%%%%%%%%%%%%%%%%%%%%%%%%%%%%%%%%%%%%
\begin{Example}\label{ex:colorglnv} Let $A=\bigoplus_{\al\in G}A_\al $ be a $G$-graded associative algebra with
multiplication
$A_\al A_\be\subseteq A_{\al+\be}$.
Define the bracket
$$[x, y]:=xy- \ve(x,y)y x,$$
for all $x\in A_{\al},y\in A_{\be}$. Then $(A, [\cdot,\cdot])$ is a Lie color algebra.
\end{Example}

Let $G$ be an abelian group and $V = \bigoplus_{\al\in G} V_\al$  a $G$-graded vector
space. Then the associative algebra $\End_G V$ is equipped with the induced $G$-grading
$\End_G (V) =\bigoplus_{\al\in G} \End_\al V$, where
$$\End_\al (V)  = \{A\in \End (V) |A(V_\al)\subseteq V_{\gamma+\al}\}.$$
By Example \ref{ex:colorglnv} we get a Lie color algebra on $(\End_G(V), [\cdot,\cdot])$,
which is denoted by $\glnv$.

A homomorphism  between two Lie color algebras $(\g,[\cdot,\cdot ])$ and $(\g',[\cdot,\cdot]')$ is a grade-perserving
linear map $\varphi: \g \to \g'$ such that
$$\varphi([x,y]) = [\varphi(x),
\varphi(y)]' $$ %\quad \forall x,y\in L, \,  \, \forall \al\in G.$$
for all $x, y\in  \g$.

Let $\g$ be a Lie color algebra and  $V$ is a $G$-graded vector space, a {\bf representation} or a {\bf $\g$-module} on $V$ is a homomorphism $\rho: \g\to \glnv$ such that
\begin{align}\label{colormodule}
\rho([x, y])v = \rho(x)(\rho(y)v) -\ve(x, y)\rho(y)(\rho(x)v),
\end{align}
for all $x\in \g_{\al},y\in \g_{\be},v\in V$.

Recall that, the cohomology of Lie color algebra $\g$ with coefficients in a representation $V$ is defined as follows:
The cochain complex  $C^{n}(\g,V)$ is the space of
$\ve$-skew symmetric homomorphisms from $\g^{\otimes n}$ to $V$, $n\geq 0$ and $C^{0}(\g,V)=V$.
The coboundary operator
$\delta:C^{n-1}(\g,V)\to C^{n}(\g,V)$
is defined by
\begin{eqnarray}\label{eq:Liecohom}
\notag&&\delta\omega(x_1,x_2,\cdots,x_n)\\
\notag&=&\sum_{k=1}^{n}(-1)^{k+1}\ve(x_{1}+\cdots+x_{k-1},x_k)[x_k,\omega(x_1,\cdots,\widehat{x_k},\cdots,x_n)]\\
\notag&&+\sum_{1\leq k<j\leq n}(-1)^{j+1}\ve(x_{k+1}+\cdots+x_{j-1},x_j)
\omega(x_1,\cdots,x_k\circ x_j,\cdots, \widehat{x_j},\cdots,x_n),
\end{eqnarray}
where $\omega\in C^{n-1}(\g,V)$ and $x_k\in V$.
It is proved that $\delta\circ \delta=0$, see \cite{SZ}.
Therefore $\{C^*(\g,V), \delta\}$ is indeed a cochain complex, whose cohomology is called the
cohomology of the Lie color algebra $\g$ with coefficients in the representation $V$.

\begin{Example}\label{example1}\cite{OS,CSO}
The $G = \mathbb{Z}_2^2$-graded analogue of $\mathfrak{sl}_2(\mathbb{C})$ is defined as an algebra $\g=\oplus_{\al\in \mathbb{Z}^2_2}\g_\al$ with
grading
$$\g_{(0,0)}= 0, \quad \g_{(1,0)}=\mathbb{C}x, \quad \g_{(0,1)}=\mathbb{C}y,\quad \g_{(1,1)} = \mathbb{C}z.$$
and bicharacter
$$\varepsilon((\alpha_1, \alpha_2), (\beta_1, \beta_2)) := (-1)^{\alpha_1\beta_2-\alpha_2\beta_1}$$
for all $\alpha_1, \alpha_2, \beta_1, \beta_2\in \mathbb{Z}_2$ where the basis elements $x,y,z$ satisfying the following relations
\begin{eqnarray*}
&&[x,x]=[y,y]=[z,z]=0,\\
&&[x, y]=xy+y x=z,\\
&& [y, z]=yz+z y=x,\\
&& [z, x]=zx+xz=y.
\end{eqnarray*}
We  denote this three-dimensional Lie color algebra by $\mathfrak{sl}_2^c$. The cohomology of this algebra was investigated in \cite{PS,SZ}.
\end{Example}

\section{Lie Color 2-algebras and 2-term Color $L_{\infty}$-algebras}

In this section, we introduced the concept of Lie color 2-algebras and 2-term color $L_{\infty}$-algebras as an generalization of Lie 2-algebras.
It is proved that the category of Lie color $2$-algebras and the category of $2$-term color $L_{\infty}$-algebras are equivalent.

%\subsection{Lie color 2-algebras}

Recall that we denote the category of $G$-graded vector spaces or color vector spaces by ${\rm\mathbf{Vect}}^G$.
%We also call $G$-graded vector spaces by color vector spaces.

\begin{Definition}
A {\bf color 2-vector space} is a category in ${\rm\mathbf{Vect}}^G$.
\end{Definition}

Thus, a color $2$-vector space $\calV$ is a category with a $G$-graded vector space of
objects $\calV_0=\bigoplus_{\al\in G}{(\calV_0)}_\al $ and a $G$-graded vector space of morphisms $\calV_1=\bigoplus_{\al\in G}{(\calV_1)}
_\al $, such that the
source and target maps $s,t \maps \calV_{1} \rightarrow \calV_{0}$, the
identity-assigning map $i \maps \calV_{0} \rightarrow \calV_{1}$, and the
composition map $\circ \maps \calV_{1} \times _{\calV_{0}} \calV_{1}
\rightarrow \calV_{1}$ are all grade-preserving linear maps.  As in the ordinary case in \cite{Baez}, we write a
morphism $f$ from source $x$ to target $y$ by $f \maps x \to y$, i.e. $s(f) = x$ and $t(f) = y$.
We also write $i(x)$ as $1_x$.

Color 2-vector spaces are in one-to-one correspondence with 2-term complexes of color vector spaces.
A 2-term complex of color vector spaces is a pair of color vector space with a differential between them: $\calC_1\stackrel{d}
{\longrightarrow}\calC_0$.
Roughly speaking, given a color 2-vector space $\calV$,
$\Ker(s)\stackrel{t}{\longrightarrow}\calV_0$ is a 2-term complex.
Conversely, any 2-term complex of color vector spaces
$\calV_1\stackrel{d}{\longrightarrow}\calV_0$ gives rise to a
color 2-vector space of which the set of objects is $\calC_0$, the set of
morphisms is $\calC_0\oplus \calC_1$, the source map $s$ is given by
$s(x,h)=x$, and the target map $t$ is given by $t(x,h)=x+dh$,
where $x\in \calV_0,~h\in \calV_1.$ We denote the color 2-vector space associated
to the 2-term complex of color vector spaces
$\calV_1\stackrel{d}{\longrightarrow}\calV_0$ by $\calV$:
\begin{equation}\label{eqn:V}
\calV=\begin{array}{c}
\calV_1:=\calC_0\oplus \calC_1\\
\vcenter{\rlap{s }}~\Big\downarrow\Big\downarrow\vcenter{\rlap{t }}\\
\calV_0:=\calC_0.
 \end{array}\end{equation}

\begin{Definition} \label{defnlie2alg}
A {\bf Lie color $2$-algebra} consists of a  color $2$-vector space $\calL$ equipped with
\begin{itemize}
\item a $\ve$-skew-symmetric bilinear functor, the {\bf bracket}, $[\cdot, \cdot]\maps \calL \times \calL
\rightarrow \calL$
\item a completely $\ve$-skew symmetric
trilinear natural isomorphism, the {\bf Jacobiator},
$$J_{x,y,z} \maps [x,[y,z]] \to [[x,y],z] +\ve(x,y) [y,[x,z]],$$
\end{itemize}
such that the following  {\bf Jacobiator identity} is satisfied
\begin{align}\label{Jacobiator}
&\quad J_{x,y,[z,t]} \Big(1+\ve(x,y)[y,J_{x,z,t}]\Big) \Big({J_{[x,y],z,t}+\ve(x,y)J_{y,[x,z],t}+\ve(x,y+z)J_{y,z,[x,t]}}\Big)\notag\\
&=[x,J_{y,z,t}] \Big(J_{x, [y,z],t}+\ve(y,z)J_{x,z,[y,t]}\Big) \Big([J_{x,y,z},t]+ 1+1+\ve(x+y,z)[z,J_{x,y,t}] \Big).
\end{align}
\end{Definition}
\footnote{The above Jacobiator is in fact a left Jacobiator since we can define a right one by $J'_{x,y,z} \maps [[x,y],z] \to [x,[y,z]] +\ve(y,z) [[x,z],y]$. From it, we get a right Jacobiator identity:
\begin{align}\label{Jacobiator'}
&\quad J_{[w,x],y,z} \Big(1+\ve(y,z)[J_{w,x,z},y]\Big) \Big(J_{w,x, [y,z]}+\ve(y,z)J_{w, [x,z], y} +\ve(x+y,z) J_{[w,z],x,y}\Big)\notag\\
&=[J_{w,x,y},z] \Big(J_{w, [x,y],z}+\ve(x,y)J_{[w,y],x,z}\Big) \Big([w, J_{x,y,z}]+ 1+1+\ve(x,y+z)[J_{w,y,z},x] \Big).
\end{align}
}
When we draw the Jacobiator identity as a commutative diagram, we see that it relates two
ways of using the Jacobiator to rebracket the expression
$[x,[y,[z,t]]]$:
$$\def\objectstyle{\scriptstyle}
  \def\labelstyle{\scriptstyle}
\xymatrix{
&[x,[y,[z,t]]]\ar[dr]^{J_{x,y,[z,t]}}\ar[dl]_{[x,J_{y,z,t}]}&\\
[x,[[y,z],t]]+\ve(y,z)[x,[z,[y,t]]]\ar[dd]^{J_{x,[y,z],t}+\ve(y,z)J_{x,z,[y,t]}}
&& [[x,y],[z,t]]+\ve(x,y)[y,[x,[z,t]]]\ar[dd]_{1+\ve(x,y)[y,J_{x,z,t}]}\\
&&\\
{\begin{aligned}&\scriptstyle\quad[[x,[y,z]],t]+\ve(x,y+z) [[y,z],[x,t]]\\[-.5em]
&\scriptstyle+\ve(y,z)[[x,z],[y,t]]+\ve(y,z)\ve(x,z)[z,[x,[y,t]]]\end{aligned}}
\ar[dr]_{[J_{x,y,z},t]+1+1+\ve(x+y,z)[z,J_{x,y,t}]}&&
{\begin{aligned}&\scriptstyle[[x,y],[z,t]]+\ve(x,y)[y,[[x,z],t]] \\[-.5em]
&\scriptstyle \quad+\ve(x,y)\ve(x,z) [y,[z,[x,t]]] \end{aligned}}
\ar[dl]^{J_{[x,y],z,t}+\ve(x,y)J_{y,[x,z],t}+\ve(x,y+z)J_{y,z,[x,t]}}\\
&P=Q&}
\\ \\
$$
where $P$ and $Q$ are given by
\begin{eqnarray*}
  P&=&[[[x,y],z],t]+\ve(x,y)[[y,[x,z]],t]\\
  &&+\ve(x,y+z) [[y,z],[x,t]]+\ve(y,z)[[x,z],[y,t]]\\
   &&+\ve(y,z)\ve(x,z)[z,[[x,y],t]]+\ve(y,z)\ve(x,z)\ve(x,y)[z,[y,[x,t]]],\\
  Q&=&[[[x,y],z],t]+\ve(y,z)\ve(x,z)[z,[[x,y],t]]\\
  &&+\ve(x,y)[[y,[x,z]],t]+\ve(x,y)\ve(y,x+z)[[x,z],[y,t]]\\
  &&+\ve(x,y)\ve(x,z) [[y,z],[x,t]]+\ve(x,y)\ve(x,z)\ve(y,z)[z,[y,[x,t]]].
\end{eqnarray*}
They are equal since $\ve$ is a bicharacter.

\begin{Definition}
  Given Lie color 2-algebras $(\calL,\br,, \Phi)$ and $(\calL^\prime,\br,', \Phi^\prime)$, a Lie color 2-algebra morphism $F:\calL\longrightarrow
  \calL^\prime$ consists of:
  \begin{itemize}
   \item[$\bullet$] a functor $(F_0,F_1)$ from the underlying
   color 2-vector space of $\calL$ to that of $\calL^\prime$.

 \item[$\bullet$] a $\ve$-skew-symmetric natural transformation
 $$
F_2(x,y):[F_0(x),F_0(y)]'\longrightarrow F_0([x,y])
 $$
such that the following diagram commutes:
  \end{itemize}
 $$
\footnotesize{ \xymatrix{
[F_0(x),[F_0(y),F_0(z)]']'\ar[d]_{J_{F_0(x),F_0(y),F_0(z)}}\ar[rr]^{\qquad\qquad\qquad[1,F_2(y,z)]'
\quad\qquad\quad}&&[F_0x,F_0[y,z]]'\ar[d]^{F_2(x,[y,z])}\\
~[[F_0(x),F_0(y)]',F_0(z)]'+\ve(x,y)[F_0(y),[F_0(x),F_0(z)]']'\ar[d]_{[F_2(x,y),1]'+[1,F_2(x,z)]'}&&F_0[x,[y,z]]\ar[d]^{F_1J_{x,y,z}}\\
~[F_0([x,y]),F_0(z)]'+\ve(x,y)[F_0(y),F_0[x,z]]'\ar[rr]^{\footnotesize{F_2([x,y],z)+F_2(y,[x,z])}}&&
\footnotesize{F_0[[x,y],z]+\ve(x,y)F_0[y,[x,z]]. }}}
$$
\end{Definition}

The identity morphism ${\id}_\calL:\calL\longrightarrow \calL$ has the identity
functor as its underlying functor, together with an identity natural
transformation as $({\id}_\calL)_2$.  Let $\calL,~\calL'$ and $\calL''$ be Lie color
2-algebras, the composition of a pair of Lie color 2-algebra morphisms
$F:\calL\longrightarrow \calL'$ and $G:\calL'\longrightarrow \calL''$, which we
denote by $G\circ F$, is given by letting the functor $((G\circ
F)_0,(G\circ F)_1)$ be the usual composition of $(G_0,G_1)$ and
$(F_0,F_1)$, and letting $(G\circ F)_2$ be the following composite:
$$
 \xymatrix{
[G_0\circ F_0(x),G_0\circ F_0(y)]''\ar[dd]_{G_2(F_0(x),F_0(y))}\ar[dr]^{(G\circ F)_2(x,y)}&&\\
&G_0\circ F_0[x,y].&\\
 G_0[F_0(x),F_0(y)]'\ar[ur]_{G_1(F_2(x,y))}&&,}.
$$

It is easy to see that
\begin{Proposition}
There is a category {\bf LieC2Alg} with Lie color 2-algebras as objects and Lie color 2-algebra morphisms as morphisms.
\end{Proposition}

%\subsection{2-term color $L_{\infty}$-algebras}% and Omni-Lie color Algebras}

\begin{Definition} \label{2termliealgebra}
A $2$-term color $L_{\infty}$-algebra $\calV=\calV_0\oplus \calV_1$ is a complex consisting of the following
data:
\begin{itemize}
  \item two color vector spaces $\calV_{0}$ and
   $\calV_{1}$ together with a grade-preserving linear map
   $d\maps V_{1} \rightarrow V_{0}, \ d((\calV_1)_\al )\subseteq (\calV_0)_\al $.

  \item a bilinear map $l_{2}\maps \calV_{i} \times \calV_{j}
   \rightarrow \calV_{i+j},$ where $0 \leq i + j \leq 1$,

  \item a trilinear map $l_{3}\maps \calV_{0} \times \calV_{0} \times
   \calV_{0} \rightarrow \calV_{1}.$
\end{itemize}

These maps satisfy the following conditions:
\begin{itemize}
  \item[(a)] $l_2(x,y) +\ve(x,y)l_2(y,x)=0$,
  \item[(b)] $l_2(x,h) +\ve(x,h)l_2(h,x)=0$,
  \item[(c)] $l_2(h,k)=0$,
  \item[(d)] $l_{3}(x,y,z)$ is totally $\ve$-skew symmetric,
  \item[(e)] $d(l_2(x,h)) = l_2(x,dh)$,
  \item[(f)] $l_2(dh,k) =  l_2(h,dk)$,
  \item[(g)] $d(l_{3}(x,y,z))=l_2(x,l_2(y,z))-l_2(l_2(x,y),z)-\ve(x,y)l_2(y,l_2(x,z))$,
  \item[(h)] $l_{3}(x,y,dh)=l_2(x,l_2(y,h))-l_2(l_2(x,y),h)-\ve(x,h)l_2(y,l_2(x,h))$,
  \item[(i)] $\delta l_3(x,y,z,t):=l_2(x, l_3(y, z,t)) - \ve(x,y) l_2(y, l_3(x, z,t))$
    \begin{flalign*}
     &+\ve(x+y,z)l_2(z, l_3(x, y,t))+l_2(l_3(x, y, z),t)- l_3(l_2(x, y), z,t)&\\
     &+ \ve(y,z)l_3(l_2(x, z),y,t)-\ve(x,y+z)l_3(y, z,l_2(x,t))&\\
     &+ l_3(x,l_2(y, z), t) + \ve(y,z)l_3(x,z,l_2(y,t))- l_3(x, y,l_2(z,t))=0.&
      \end{flalign*}
\end{itemize}
for all homogeneous elements $x,y,z,t\in \calV_{0}$ and $h, k \in \calV_{1}.$
\end{Definition}

%\begin{itemize}
%\item[($g'$)] $\ve(z,x)d(l_{3}(x,y,z))=\ve(z,x)[[x,y],z]+\ve(x,y)[[y,z],x]+\ve(y,z)[[z,x],y]=J_1(x,y,z)$,
%\item[($h'$)] $\ve(z,h)l_{3}(x,y,dh)=\ve(h,x)[[x,y],h]+\ve(x,y)[[y,h],x]+\ve(y,h)[[y,x],h]=J_1(x,y,h)$,
%\end{itemize}

\begin{Definition}\label{defi:Lie-2hom}
Let $(\calV;\dM,l_2,l_3)$ and $(\calV';\dM',l_2,l_3')$ be two $2$-term color $L_{\infty}$-algebras.
A $CL_\infty$-morphisms $f$ from $\calV$ to $\calV'$ consists of
 linear maps $f_0:\calV_0\rightarrow \calV_0',~f_1:\calV_{-1}\rightarrow \calV_{-1}'$
 and $f_{2}: \calV_{0}\times \calV_0\rightarrow \calV_{-1}'$,
such that the following equalities hold for all $ x,y,z\in \calV_{0},
a\in \calV_{-1},$
\begin{itemize}
\item [$\rm(i)$] $f_0\dM=\dM'f_1$,
\item[$\rm(ii)$] $f_{0}l_2(x,y)-l'_2(f_{0}(x),f_{0}(y))=\dM'f_{2}(x,y),$
\item[$\rm(iii)$] $f_{1}l_2(x,a)-l'_2f_{0}(x),f_{1}(a))=f_{2}(x,\dM a)$,
%\item[$\rm(iii)_2$] $f_{1}[a,x]-[f_{1}(a),f_{0}(x)]'=f_{2}(\dM a,x)$,
\item[$\rm(iv)$] $f_1(l_3(x,y,z))-l_3'(f_0(x),f_0(y),f_0(z))$
  $=f_2(x, l_2(y,z))- f_2(l_2(x,y),z) - f_2(y,l_2(x,z))$ $+ l'_2(f_0(x), f_2(y,z)) - l'_2(f_2(x,y), f_0(z))- l'_2(f_0(y), f_2(x,z)).$
\end{itemize}
 If $f_2=0$, the $CL_\infty$-morphisms $f$ is called a strict $CL_\infty$-morphisms.
\end{Definition}

Let $f:\calV\to \calV'$ and $g:\calV'\to \calV''$ be two $CL_\infty$-morphisms, then their composition $g\circ f:\calV\to \calV''$
is a $CL_\infty$-morphism defined as $(g\circ f)_0=g_0\circ f_0$, $(g\circ f)_1=g_1\circ f_1$
and
$$(g\circ f)_2(x,y)=g_2 (f_0(x), f_0(y))+g_1(f_2(x,y)).$$

The identity $CL_\infty$-morphism $1_\calV: \calV\to \calV$ has the identity chain map together with $(1_\calV)_2=0$.

\begin{Proposition}
There is a category {\bf 2CL$_\infty$} with 2-term color $L_\infty$-algebras as objects and $CL_\infty$-morphisms as morphisms.
\end{Proposition}

Now we establish the equivalence between the category of Lie color $2$-algebras and $2$-term color $L_{\infty}$-algebras.

\begin{Theorem} The categories  ${\bf2CL_\infty}$ and  ${\bf LieC2Alg}$ are equivalent.
\end{Theorem}

Since the proof is similar as in \cite{Baez}, we omit the details.

\section{Construction of Lie Color 2-algebras}

In this section, special cases and concrete examples of Lie color 2-algebras are given.
This include Lie color algebras with  3-cocyles, crossed module of  Lie color algebras, string Lie color algebra, and omni-Lie color algebras.
At last, we prove that from any Leibniz color algebra, we can obtain Lie color 2-algebras.

\subsection{Skeletal Lie color 2-algebras}

A 2-term color $L_{\infty}$-algebra is called {\bf skeletal} if $d=0$.
In this case,
from conditions $(a)$ and $(g)$, we have $\calV_0$ is a Lie color algebra.
Conditions $(b)$ and $(h)$ imply that $\calV_1$ is a representation of $\calV_0$ by the action
defined by $\rho(x)h := l_2(x, h)$. Now condition $(i)$ can be described in terms of a 3-cocycle condition in the
Lie color algebra cohomology of $\calV_0$ with values in $\calV_1$.

\begin{Proposition}
Skeletal 2-term color $L_{\infty}$-algebras are in one-to-one correspondence with
quadruples $(\g, V, \rho, l_3)$ where $\g$ is a Lie color algebra, $V$ is a $G$-graded vector space,
$\rho$ is a representation of $\g$ on $V$ and $l_3$ is a 3-cocycle on $\g$ with values in $V$.
\end{Proposition}

Recall that a quadratic Lie color algebra is Lie color algebra  $(\g,[\cdot,\cdot])$ together with an $\ve$-symmetric,
nondegenerate, invariant bilinear form $B:\g\times\g\longrightarrow\R$, such that for any $x,y,z\in\g$,
\begin{eqnarray*}
B(x,y)=\ve(x,y)B(y,x),\quad B([x,y],z)=B(x,[y,z]).
\end{eqnarray*}

\begin{Example}
Given a quadratic Lie color algebra $(\g,[\cdot,\cdot], B)$, we construct a 2-term color $L_{\infty}$-algebra as follows.
Let $V_1=\R, V_0=\g, d=0$, and define $l_2, l_3$ by
\begin{equation}\label{eqn:l2l3string}
l_2(x,y)=[x,y],\quad l_2(x,h)=0,\quad l_3(x,y,z)=B([x,y],z),
\end{equation}
where $x,y,z\in\g, h\in \R$.
All the conditions in the Definition \ref{2termliealgebra}
are satisfied, and we get a 2-term color $L_{\infty}$-algebra $(\R\stackrel{0}{\longrightarrow}\g,l_2,l_3)$
from a quadratic Lie color algebra $(\g,[\cdot,\cdot], B)$. We call this {\bf string Lie color 2-algebra}.
\end{Example}

\subsection{Strict Lie color 2-algebras}
Another kind of  2-term color $L_{\infty}$-algebra is called {\bf strict} if $l_3=0$.
This kind of Lie color 2-algebras can be described in terms of crossed modules of Lie color algebras.

\begin{Definition} Let $(\g,[\cdot,\cdot]_{\g})$ and $(\h,[\cdot,\cdot]_{\h})$ be two Lie color algebras.
A crossed module of Lie color algebras is a homomorphism of Lie color algebras
$\varphi: \h\to \g$ together with a representation of $\g$ on $\h$, denoted by $x\triangleright h:=\rho(x)h$,
such that
$$\varphi(x\triangleright h) = [x, \varphi(h)]_{\g},\quad\varphi(h)\triangleright k = [h, k]_{\h},$$
for all $h, k\in\h, x\in\g$.
\end{Definition}

\begin{Proposition}
Strict 2-term color $L_{\infty}$-algebras are in one-to-one correspondence with crossed modules of Lie color algebras.
\end{Proposition}

\pf Let $\calV_1\stackrel{d}{\longrightarrow} \calV_0$ be a $2$-term color $L_{\infty}$-algebra with $l_3=0$. We construct Lie color algebras on $\g=\calV_0$ and $\h=\calV_1$ as follows.
The bracket  on $\g$ and $\h$ are defined by
\begin{eqnarray*}
&&[h,k]_{\h}:=l_2(dh,k),\quad\forall~x,y\in \h=\calV_1;\\
&&[x,y]_{\g}:=l_2(x,y),\quad\forall~h,k\in \g=\calV_0.
\end{eqnarray*}
By condition $(a)$ and $(g)$ in Definition \ref{2termliealgebra},
it is easy to see that $[\cdot,\cdot]_\g$ satisfies the $\ve$-Jacobi identity.
By $(h)$, we have
\begin{eqnarray*}
 && [h,[k,l]_\h]_\h-[[h,k]_\h,l]_\h+\ve(k,l)[[h,l]_\h k]_\h\\
 &=&l_2(dh,l_2(dk,l))-l_2(dl_2(dh,k),l)+\ve(k,l)l_2(dl_2(dh,l),k)\\
 &=&l_2(dh,l_2(dk,l))-l_2(l_2(dh,dk),l)+\ve(k,l)l_2(l_2(dh,dl),k)\\
 &=&0.
\end{eqnarray*}
Thus $[\cdot,\cdot]_\h$ satisfies the $\ve$-Jacobi identity.
Now let $\varphi=d$, then by $(e)$, we have
$$
\varphi([h,k]_\h)=d(l_2(dh,k))=l_2(dh,dk)=[\varphi(h),\varphi(k)]_\g,
$$
which implies that $\varphi$ is a homomorphism of Lie color algebras.

Now define the maps of $\trr:\g\times \h\to \h$ by
$$x\trr h:=l_2(x,h)\in\h,$$
which is a representation of $\g$ on $\h$ and it is easy to check that
\begin{eqnarray*}
&&\varphi(x\trr h)=d(l_2(x,h)) = l_2(x,dh)= [x, \varphi(h)]_{\g}\\
&&\varphi(h)\trr k =[dh,k]=[h, k]_{\h}.
\end{eqnarray*}
Therefore, we obtain a crossed module of Lie color algebras.

Conversely,  a crossed module of Lie  color algebras
 gives rise to a 2-term color $L_{\infty}$-algebra with $d=\varphi$,
$\calV_0=\g$ and $\calV_1=\h$, where the brackets are given by
\begin{eqnarray*}
~ l_2(x,y)&:=&[x,y]_{\g},\quad \forall
~x,y\in\g;\\
~l_2(x,h)&:=&x\trr h,\quad\forall~ x\in\h;\\
~l_2(h,k)&:=&0.
\end{eqnarray*}
The crossed module conditions give various conditions for $2$-term color $L_\infty$-algebras with $l_3=0$.

\qed

\medskip

%\section{Derivations and Automorphism Groups of Lie Color Algebras}
Let $\g$ be a Lie color algebra, recall that a map $D\in\gl(\g)$ is called
a homogeneous color derivation of degree $|D|$ if $D(\g_\al )\subseteq V_{\al+|D|}$ and
$$D([x,y])=[Dx,y]+\ve(D,x)[x,Dy].$$
Denote by $\Der(\g)=\bigoplus_{\al\in G}\Der_{\al}(\g)$, where
$\Der_{\al}(\g)$ is the vector space spanned by all homogeneous color derivation of degree $\al$.
We find that $\Der(\g)$ becomes a Lie color algebra under the bracket
$$[D,D']=DD'-\ve(D,D')D'D,$$
where $D,D'$ are homogeneous color derivations of degree $|D|, |D'|$.

For a bilinear operation $\om$ on $V$ such that $\om: V_\al\times V_\be\to V_{\al+\be}$, we define the adjoint operator
$\ad_\om:V_\al\rightarrow \gl(V)_\al$ by
$$\ad_{\om}(x)(y)=\om(x,y)\in V_{\al+\be}, \forall x\in V_\al,\, y\in V_\be.$$
Then the graph of the adjoint operator
$$\calf_\om :=\{\ad_\om x+x ~; \, \forall x\in V\}\subset \cale$$
is a subspace of $\cale$.

\begin{Proposition}\label{prop:derivation}
$D$ is a  homogeneous color derivation of $V$ if and only if $\calf_\om$ is an invariant subspace of $D$, that is $D\circ\calf_\om\subseteq\calf_\om$.
\end{Proposition}
\pf For $\ad_{\om}(x)+x\in \calf_\om$,
\begin{align*}
D\circ (\ad_{\om}(x)+x)&=[D, \ad_{\om}(x)]+Dx.
\end{align*}
The right hand side is belong to $\calf_\om$ if and only if
\begin{align*}
[D,\ad_{\om}(x)]&=\ad_{\om}(Dx),
\end{align*}
that is
\begin{eqnarray*}
&& D\ad_{\om}(x)(y)-\ve(D,x) \ad_{\om}(x)D(y)-\ad_{\om}(Dx)(y)\\
&=&D[x,y]-\ve(D,x)[x, D(y)]-[Dx,y]\\
&=&0.
\end{eqnarray*}
Thus $D$ is a derivation if and only if $D\circ\calf_\om\subseteq\calf_\om$.
\qed

We call the set of elements $D\in \gl(\g)$ such that $D\circ \calf_\om\subseteq \calf_\om$ the normalizer of $\calf_\om$, which is denoted by $N(\calf_\om)$.
\begin{Proposition}\label{prop:derivation}
Let $D,D'$ be two homogeneous color derivations. Then $[D,D']$ is also a homogeneous color derivation.
Thus we have $\Der(\g)=N(\calf_\om)$ is a Lie color subalgebra of $\gl(\g)$.
\end{Proposition}
\pf Let $D\circ\calf_\om\subseteq\calf_\om$ and $D'\circ\calf_\om\subseteq\calf_\om$, then
\begin{align*}
[D,\ad_{\om}(x)]=\ad_{\om}(Dx),\quad [D',\ad_{\om}(x)]=\ad_{\om}(D'x).
\end{align*}
By the $\ve$-Jacobi identity of $\g$, we have
\begin{eqnarray*}
&&[[D,D']\circ\ad_{\om}(x)]\\
&=&[D,[D',\ad_{\om}(x)]]+\ve(D',x)[[D,\ad_{\om}(x)],D']\\
&=&[D,\ad_{\om}(D'x)]+\ve(D',x)\ve(D+x,D')[D',\ad_{\om}(Dx)]\\
&=&\ad_{\om}(DD'x)+\ve(D,D')\ad_{\om}(D'Dx)\\
&=&\ad_{\om}([D,D']x).
\end{eqnarray*}
This is equivalent to $[D,D']\circ\calf_\om\subseteq\calf_\om$, so the bracket in $\Der(\g)$ is closed.
Thus $\Der(\g)$ is a Lie color algebra as a Lie color subalgebra of $\gl(\g)$.
\qed

\begin{Example}
Let $\g$ be a Lie color algebra, $\Der(\g)$ and  $\Inn(\g)$ be the set of their derivations and inner derivations.
Then we obtain a crossed module $i:\Inn(\g)\to\Der(\g)$, with $\Der(\g)$ acting  $\Inn(\g)$ by $D\trr\ad_x=\ad_{Dx}$.
\end{Example}

\subsection{Omni-Lie color algebras}
The notion of omni-Lie algebras was generalized to omni-Lie superalgebras in \cite{ZL}.
In thus subsection, we introduce the concept of omni-Lie color algebras and construct
2-term color $L_{\infty}$-algebras from them.

%Here we give the color analogue of omni-Lie algebras, which is called omni-Lie color algebras.

Let $\g$ be a Lie color algebra and  $V$ be a $\g$-module. We define an operation $\circ$ on $\cale:=\g\oplus V$
by
\begin{equation}
  \label{eq:Leibnizbracket}
   (A+x)\circ (B+y) = [A,B]+Ay,
\end{equation}
for all $A, B\in\g$ and $x,y\in V$.
Then it is easy to check:
\begin{Proposition}
$(\cale,\circ)$ is a Leibniz color algebra.
\end{Proposition}

\emptycomment{
\pf We have to check the $\ve$-Leibniz rule for the operation $\circ$.
Let $e_1=A+x, e_2=B+y, e_3=C+z$. Then we have
\begin{eqnarray*}
&&\{e_1\circ e_2\} \circ e_3-e_1\circ \{ e_2 \circ e_3\} -\ve(x,y)e_2\circ \{ e_1 \circ e_3\} \\
&=&([A,B]+Ay)\circ (C+z)-(A+x)\circ ([B,C]+Bz)\\
&&-\ve(x,y)(B+y)\circ ([A,C]+Az)\\
&=&[[A,B],C]-[A,[B,C]]-\ve(x,y)[B,[A,C]]\\
&&+[A,B]z-ABz-\ve(x,y)BAz\\
&=&0.
\end{eqnarray*}
The right hand side of the above equation is zero because $\g$ is a Lie color algebra acting on $V$.
\qed
}

We call this type of Leibniz color algebra structure on $\cale$  the hemisemidirect product of $\g$ with $V$ as in \cite{KW},
denote it by $\g\ltimes_H V$.
Note that the operation $\circ$ is not $\ve$-skew-symmetry.
The $\ve$-skew-symmetrized bracket of it is:
\begin{equation}  \label{eq:bracket}
  \bleft A+x, B+y\bright := [A,B]+\half\left(Ay-\ve(x,y)Bx\right).
\end{equation}
This is called demisemidirect product of $\g$ with $V$, denoted by $\g\ltimes_D V$.
Furthermore, we define a $\ve$-symmetric bilinear
form on $\cale$ with values in $V$ by
\begin{equation}  \label{eq:symmetric}
  \langle A+x, B+y\rangle:=\half(Ay+\ve(x,y)Bx).
\end{equation}
The triple $(\cale, \bleft\cdot,\cdot\bright, \la\cdot,\cdot\ra)$ is called an {\bf omni-Lie  color algebra}.

\begin{Proposition}\label{prop:homotopy}
For $e_1=A+x, e_2=B+y, e_3=C+z\in \cale$, define
\begin{align*}
J_1(e_1, e_2, e_3):=&\ \ \ve(z,x) \bleft\bleft e_1,e_2\bright, e_3\bright+ \ve(x,y)\bleft\bleft
e_2,e_3\bright, e_1\bright+ \ve(y,z)\bleft\bleft e_3,e_1\bright, e_2\bright,\\
T(e_1, e_2, e_3):=&\three\Big(\ve(z,x) \langle\bleft e_1,e_2\bright, e_3\rangle+ \ve(x,y)\langle\bleft
e_2,e_3\bright, e_1\rangle+ \ve(y,z)\langle\bleft e_3,e_1\bright, e_2\rangle\Big),\\
S(e_1, e_2, e_3):=&\four\Big(\ve(z,x)(e_1\circ e_2)\circ e_3+\ve(x,y)(e_2\circ e_3)\circ e_1+\ve(y,z)(e_3\circ e_1)\circ e_2\Big).
\end{align*}
Then we have
$$J_1(e_1, e_2, e_3)= T(e_1, e_2, e_3)=S(e_1, e_2, e_3).$$
\end{Proposition}
\pf  We compute $J_1$ and $T$ as follows:
\begin{eqnarray*}
&&J_1(e_1, e_2, e_3)\\
&=&\ve(z,x)\bleft\bleft A+x, B+y\bright, C+z\bright+\mbox{c.p.}\\
&=&\bleft\ve(z,x)[A,B]+\half\ve(z,x)\left(Ay-\ve(x,y)Bx\right), C+z\bright+\mbox{c.p.}\\
&=&\ve(z,x)[[A,B],C]+\ve(x,y)[[B,C],A]+\ve(y,z)[[C,A],B]\\
&&+\half\left(\ve(z,x)[A,B]z-\half\ve(z,x)\ve(x+y,z)C\left(Ay-\ve(x,y)Bx\right)\right)\\
&&+\half\left(\ve(x,y)[B,C]x-\half\ve(x,y)\ve(y+z,x)A\left(Bz-\ve(y,z)Cy\right)\right)\\
&&+\half\left(\ve(y,z)[C,A]y-\half\ve(z,y)\ve(z+x,y)B\left(Cx-\ve(z,x)Az\right)\right)\\
&=&\four\Big(\ve(z,x)[A,B]z+\ve(x,y)[B,C]x+\ve(y,z)[C,A]y\Big).
\end{eqnarray*}
and
\begin{eqnarray*}
  &&T(e_1, e_2, e_3)\\
  &=&\textstyle{\frac{1}{3}}\ve(z,x)\langle\bleft A+x,B+y\bright, C+z\rangle +\mbox{c.p.}\\
  &=&\textstyle{\frac{1}{3}}\ve(z,x)\langle[A,B]+\frac{1}{2}\left(Ay-\ve(x,y)Bx\right), C+z\rangle +\mbox{c.p.}\\
  &=&\six\Big(\ve(z,x)[A,B]z+\frac{1}{2}\ve(y,z)C\left(Ay-\ve(x,y)Bx\right)\Big) +\mbox{c.p.}\\
  &=&\six\Big(\ve(z,x)[A,B]z+\ve(x,y)[B,C]x+\ve(y,z)[C,A]y\Big)\\
  &&+\twelve\Big(\ve(z,x)[A,B]z+\ve(x,y)[B,C]x+\ve(y,z)[C,A]y\Big)\\
  &=&\four\Big(\ve(z,x)[A,B]z+\ve(x,y)[B,C]x+\ve(y,z)[C,A]y\Big).
\end{eqnarray*}
Thus, both of them are equal to $S(e_1, e_2, e_3)$ as desired.
\qed
\medskip

Now, for a Lie color algebra $\g$ and a $\g$-module $V$, let
$$\calV_0=\g\ltimes_D V,\quad  \calV_1=V,\quad d=i: V\hookrightarrow \g\ltimes_D V$$
where $i$ is the inclusion map  and define
\begin{align*}
l_2=\bleft\cdot,\cdot\bright,\quad l_3=-\ve(x,z)J_1.
\end{align*}

\begin{Theorem}\label{omniare2term}
With notations above, from an omni-Lie color algebra  $(\cale,
\bleft\cdot,\cdot\bright, \la\cdot,\cdot\ra)$ we obtain a 2-term color $L_{\infty}$-algebra $(V\stackrel{i}{\hookrightarrow}\g\ltimes_D V,\, l_2, \,
l_3)$.
\end{Theorem}

\pf
It can be checked that various conditions in Definition \ref{2termliealgebra} hold.
For example, by the grading in $\g\op V$
we have $\deg(A+x)=\deg (A)=\deg (x)$, then we have
\begin{eqnarray*}
&&\bleft A+x, B+y\bright+\ve(x, y)\bleft B+y, A+x\bright\\
&=& [A,B]+\half\left(Ay-\ve(x,y)Bx\right)+\ve(A,B)[B,A]\\
&&+\ve(x, y)\half\left(Bx-\ve(y,x)Ay\right)\\
&= &[A,B]+\ve(x,y)[B,A]+\half\left(Ay-\ve(x,y)Bx\right)\\
&&+\half\left(\ve(x, y)Bx-Ay\right)\\
&=&0.
\end{eqnarray*}
Thus condition $(a)$ holds.

Let $e_1=A, e_2=B, e_3=C$ where $A, B, C\in \g$ and $e_4=t\in V$, then we have
\begin{eqnarray*}
%&&\delta l_3(A,B,C,t)\\
&&\bleft A, l_3(B, C,t)\bright - \ve(A,B) \bleft B, l_3(A, C,t)\bright \\
&& +\ve(A+B,C)\bleft C, l_3(A, B,t)\bright+\bleft l_3(A, B, C),t\bright- l_3(\bleft A, B\bright, C,t)\\
&&+ \ve(B,C)l_3(\bleft A, C\bright,B,t)-\ve(A,B+C)l_3(B, C,\bleft A,t\bright)\\
&&+ l_3(A,\bleft B, C\bright, t) + \ve(B,C)l_3(A,C,\bleft B,t\bright)- l_3(A, B,\bleft C,t\bright)\\
%%%%%%%%%%%%%%%%%%%%%-%%%%%%%%%%%%%%%%%%%%%%-%%%%%%%%%%%%%%%%%%%%%%
&=&-\eight A[B,C]t +\eight\ve(A,B)B[A,C]t-\eight\ve(A+B,C)C[A, B]t+0\\
&& + \four[[A, B],C]t- \four\ve(B,C)[[A, C],B]t+ \four \ve(A,B)\ve(A,C)[[B,C],A]t \\
&&+\eight\ve(A,B+C)[B,C]At-\eight \ve(B,C)[A,C]Bt+\eight [A,B]Ct\\
&=&-\textstyle{\frac{3}{8}}\Big([A,[B,C]] - [[A,B],C] - \ve(A,B)[B,[A,C]]\Big)t\\
&=&0.
\end{eqnarray*}
Thus condition $(i)$ holds.
The other conditions can be checked similarly.
\qed

\begin{Example}\label{example2}
The Lie superalgebra $\mathfrak{gl}(1|1)$ consists of  $\g_0 = \mathbb{C}h_1 \oplus \mathbb{C}h_2$  and $\g_1 = \mathbb{C}e \oplus \mathbb{C}f$
where the basis elements is given by
$$h_1=\left(\begin{array}{cc} 1 & 0 \\ 0 & 0 \\ \end{array} \right),\quad
h_2=\left(\begin{array}{cc} 0 & 0 \\ 0 & 1 \\ \end{array} \right),\quad e=\left(\begin{array}{cc} 0 & 1 \\ 0 & 0 \\ \end{array} \right), \quad
f=\left(\begin{array}{cc} 0 & 0 \\ 1 & 0 \\ \end{array} \right).$$
and the bracket is given by
\begin{eqnarray*}
&&[e, f]=ef+fe=h_1+h_2,\\
&& [h_1, e]=e, \quad [h_2, e]=-e,\\
&& [h_1, f]=-f, \quad [h_2, f]=f.
\end{eqnarray*}
Let $V=V_0\oplus V_1=\mathbb{C}x_0\oplus \mathbb{C}x_1$ be the two-dimensional representation:
$$x_0=\left(
       \begin{array}{c}
         1 \\
         0 \\
       \end{array}
     \right)
, \quad x_1=\left(
       \begin{array}{c}
         0 \\
         1 \\
       \end{array}
     \right)$$
where the action of $\mathfrak{gl}(1|1)$ is given by matrix multiplication:
\begin{eqnarray*}
ex_0  = 0, \quad ex_1 = x_0 , \quad  fx_0  = x_1 , \quad fx_1  = 0,\\
h_1x_0  = x_0 , \quad h_1x_1  = 0,\quad  h_2x_0  = 0, \quad h_2x_1  = x_1 .
\end{eqnarray*}
Thus we get a 2-term super $L_{\infty}$-algebra as follows:
$$(V\stackrel{i}{\hookrightarrow}\mathfrak{gl}(1|1)\ltimes_D V,\, l_2, \,l_3)$$
where $l_2$ is equal to the bracket defined as above on $\mathfrak{gl}(1|1)$, $l_2$ is zero on $V$, the other case is given by
\begin{eqnarray*}
&&l_2(e, x_0)=0,\quad l_2(e, x_1)=\half x_0,\\
&&l_2(f, x_0)= x_1, \quad l_2(f, x_0)=0,\\
&&l_2(h_1, x_0)=\half x_0, \quad l_2(h_2, x_1)=0\\
&&l_2(h_2, x_0)=0, \quad l_2(h_2, x_1)=\half x_1.
\end{eqnarray*}
By direct computations, $l_3$ is given by
\begin{eqnarray*}
&&l_3(e,   f, x_0)=-\four x_0,  \quad l_3(e,   f, x_1)=-\four x_1,\\
&&l_3(h_1, e, x_0)= 0, \quad l_3(h_1, e, x_1)=-\four x_0,\\
&&l_3(h_2, e, x_0)= 0, \quad l_3(h_2, e, x_1)=\four x_0,\\
&&l_3(h_1, f, x_0)= \four x_1, \quad l_3(h_1, f, x_1)=0,\\
&&l_3(h_2, f, x_0)= -\four x_1, \quad l_3(h_2, f, x_1)=0.
\end{eqnarray*}
\end{Example}

\subsection{Skew-symmetrization of Leibniz color algebras}
From above subsection, we have seen that an omni-Lie color algebra $\cale=\g\ltimes_D V$
is in fact the skew-symmetrization of Leibniz color algebra $\g\ltimes_H V$ and $V$ is the left center of $\cale$.
Motivated by this, we may ask, from any Leibniz color algebra, can we get a Lie color 2-algebra? In this subsection, we give a positive answer to this question.
First we recall some concept on Leibniz color algebras.

Let $L$ be a Leibniz color algebra. We define the Leibniz kernel $\Ker(L)$ of $L$ to be the set of elements spanned by
$$\{x\circ y+\ve(x,y) y\circ x| \forall x, y\in L\}.$$
Note that if $L$ is a Lie color algebra, then $\Ker(L)=0$, otherwise if $L$ is a non-Lie Leibniz color algebra, then $\Ker(L)\neq 0$.

The left center  $Z^l(L)$ of $L$ is defined by
$$Z^l(L) = \{t\in L|t \circ x = 0, \forall x\in L)\}.$$
It is easy to see that $\Ker(L)$ and $Z^l(L)$ are ideals of $L$ and the quotient algebras $L/\Ker(L)$ and $L/Z^l(L)$ are Lie color algebras.
\begin{Proposition}\label{pro:J0}
The Leibniz kernel $\Ker(L)$ is contained in the left center $Z^l(L)$.
Thus for any non-Lie Leibniz color algebra, the set $Z^l(L)$ is not empty.
\end{Proposition}
\pf Let $x \circ y + \ve(x,y)y \circ x\in \Ker(L)$, then
\begin{eqnarray*}
&&\Big(x \circ y + \ve(x,y)y \circ x\Big) \circ z \\
&=& x \circ (y \circ z)- \ve(x,y)y \circ (x \circ z) \\
&&+\ve(x,y) \Big(y \circ (x \circ z) - \ve(y,x)  x \circ (y \circ z)\Big) \\
&=& x \circ (y \circ z)- \ve(x,y)y \circ (x \circ z) \\
&&+\ve(x,y)y \circ (x \circ z) -  x \circ (y \circ z) \\
&=& 0,
\end{eqnarray*}
for all $z\in L$. Thus $x \circ y + \ve(x,y)y \circ x\in Z^l(L)$. Therefore $\Ker(L)$ is contained in $Z^l(L)$.
\qed
\medskip

For a Leibniz color algebra  $(L,\circ)$, since the operation $\circ$ is not $\ve$-skew-symmetry, we introduce the
following $\ve$-skew-symmetric bracket on $L$ by
\begin{equation}
\Dorfman{x,y}=\half\left(x\circ y-\ve(x,y)y\circ x\right),\quad\forall x,y\in L,
\end{equation}
and let $J_{1}, S$ be given by
\begin{align}
  J_{1}(x,y,z):=&\ \ \ve(z,x)\Dorfman{\Dorfman{x,y},z}+\ve(x,y)\Dorfman{\Dorfman{y,z},x}+\ve(y,z)\Dorfman{\Dorfman{z,x},y},\\
  S(x,y,z):=&\four\Big(\ve(z,x)(x\circ y)\circ z+\ve(x,y)(y\circ z)\circ x+\ve(y,z)(z\circ x)\circ y\Big).
\end{align}

\begin{Proposition}\label{pro:J}
  Let $(L,\circ)$ be a Leibniz color algebra. Then we have
  \begin{equation}\label{eq:Leibniz}
   J_1(x,y,z)=S(x,y,z).
  \end{equation}
 \end{Proposition}
 \pf The proof is by direct computations:
\begin{eqnarray*}
&&J_1(e_1, e_2, e_3)\\
&=&\ve(z,x)\four\Big((x\circ y)\circ z-\ve(x,y)(y\circ x\big)\circ z\\
&&-\ve(x+y,z)z\circ\big(x\circ y)+\ve(x+y,z)\ve(x,y)z\circ(y\circ x\big)\Big)+\mbox{c.p.}\\
&=&\four\Big(\ve(z,x)(x\circ y)\circ z\underline{-\ve(z,x)\ve(x,y)(y\circ x)\circ z}\\
&&\underbrace{-\ve(y,z)z\circ(x\circ y)}+\underline{\underline{\ve(y,z)\ve(x,y)z\circ(y\circ x)}}\Big)\\
&&+\four\Big(\ve(x,y)(y\circ z)\circ x\underline{\underline{-\ve(x,y)\ve(y,z)(z\circ y)\circ x}}\\
&&\underline{-\ve(z,x)x\circ(y\circ z)}+\underbrace{\ve(z,x)\ve(y,z)x\circ(z\circ y)}\Big)\\
&&+\four\Big(\ve(y,z)(z\circ x)\circ y\underbrace{-\ve(y,z)\ve(z,x)(x\circ z)\circ y}\\
&&\underline{\underline{-\ve(x,y)y\circ(z\circ x)}}+\underline{\ve(x,y)\ve(z,x)y\circ(x\circ z)}\Big)\\
&=&\four\Big(\ve(z,x)(x\circ y)\circ z +\ve(x,y)(y\circ z)\circ x+ \ve(y,z)(z\circ x)\circ y\Big)\\
&=&S(x,y,z).
\end{eqnarray*}
The underline terms are canceled out, thus we obtain the result.
\qed

\begin{Proposition}\label{pro:J2}
Let $(L,\circ)$ be a Leibniz color algebra.
Then $J_{1}(x,y,z)=S(x,y,z)$ is contained in the left center $Z^l(L)$ for all $x,y,z\in L$.
 \end{Proposition}

\pf
Let $w\in L$, then by Proposition \ref{pro:J0} and Proposition \ref{pro:J} we have
\begin{eqnarray*}
&&J_{1}(x,y,z)\circ w=S(x,y,z)\circ w\\
 &=&\four\Big(\ve(z,x)(x\circ y)\circ z+\ve(x,y)(y\circ z)\circ x+\ve(y,z)(z\circ x)\circ y\Big)\circ w\\
 &=&\four\Big(\ve(z,x)x\circ (y\circ z)-\ve(z,x)\ve(x,y)y\circ (x\circ z)\\
 &&+\ve(x,y)(y\circ z)\circ x+\ve(y,z)(z\circ x)\circ y\Big)\circ w\\
  &=&\four\Big(\ve(z,x)x\circ (y\circ z)+\ve(x,y)y\circ (z\circ x)\\
 &&+\ve(z,x)\ve(x,y+z)(y\circ z)\circ x+\ve(x,y)\ve(y,z+x)(z\circ x)\circ y\Big)\circ w\\
   &=&\four\ve(z,x)\Big(x\circ (y\circ z)+\ve(x,y+z)(y\circ z)\circ x\Big)\circ w\\
 &&+\four\ve(x,y)\Big(y\circ (z\circ x)+\ve(y,z+x)(z\circ x)\circ y\Big)\circ w\\
 &=&0.
\end{eqnarray*}
Thus $J_{1}(x,y,z)\in \mathrm Z^l(L)$.
\qed
\medskip

Now it is easy to prove that
\begin{Proposition}
If $t\in \mathrm Z^l(L)$, i.e. $t\circ x=0$ for all $x\in L$, then we have
\begin{eqnarray*}
\Dorfman{x,t}=\half x\circ t,\quad S(x,y,t)=-\four\ve(t,x)(x\circ y)\circ t.
\end{eqnarray*}
\end{Proposition}

At last, for any non-Lie Leibniz color algebra $L$, we construct nontrivial Lie color 2-algebras as follows. Let
$$\calV_0=L,\quad \calV_1=\mathrm Z^l(L),\quad d=i: \mathrm Z^l(L)\hookrightarrow L,\quad l_2=\bleft\cdot,\cdot\bright,\quad l_3=-\ve(x,z)J_1.$$

\begin{Theorem}\label{thm:main1}
  With the above notations, from a Leibniz color algebra  $(L,\circ)$ we obtain a Lie color 2-algebra
   $(Z^l(L)\stackrel{d}{\hookrightarrow}L,\, l_2, \,l_3)$.
\end{Theorem}

\pf By definition of $d,~l_2$ and $l_3$ and Proposition \ref{pro:J} and Proposition \ref{pro:J2}, it is easy to see that conditions $(a)$--$(h)$ hold.
For condition $(i)$, we  verify the case of $x,y,z\in L$ and $t\in Z^l(L)$ as follows:
\emptycomment{
\begin{eqnarray*}
&& l_3(x, y,[z,t])+ \ve(x,y) [y, l_3(x, z,t)] \\
&&+ l_3([x, y], z,t)+ \ve(x,y)l_3(y, [x,z],t) +\ve(x,y+z)l_3(y, z, [x,t]) \\
&=&[x, l_3(y, z, t)]  +l_3(x,[y, z], t) +\ve(y,z)l_3(x, z,[y,t])\\
&&+[l_3(x, y, z),t]+\ve(x+y,z)[z, l_3(x, y,t)].
\end{eqnarray*}
\begin{eqnarray*}
&&[x, l_3(y, z, t)]  +l_3(x,[y, z], t) +\ve(y,z)l_3(x, z,[y,t])\\
&&+[l_3(x, y, z),t]+\ve(x+y,z)[z, l_3(x, y,t)]\\
&&- l_3(x, y,[z,t])- \ve(x,y) [y, l_3(x, z,t)] \\
&&- l_3([x, y], z,t)- \ve(x,y)l_3(y, [x,z],t) -\ve(x,y+z)l_3(y, z, [x,t]) \\
\end{eqnarray*}
}
\begin{eqnarray*}
   &&\bleft x, l_3(y, z, t)\bright - \ve(x,y) \bleft y, l_3(x, z,t)\bright\\
   &&+\ve(x+y,z)\bleft z, l_3(x, y,t)\bright+\bleft l_3(x, y, z),t\bright- l_3(\bleft x, y\bright, z,t)\\
   &&+ \ve(y,z)l_3(\bleft x, z\bright,y,t)-\ve(x,y+z)l_3(y, z,\bleft x,t\bright)\\
   && +l_3(x,\bleft y, z\bright, t) + \ve(y,z)l_3(x, z,\bleft y,t\bright)- l_3(x, y,\bleft z,t\bright)\\
   &=&-\four\Big(\bleft x, (y\circ z)\circ t\bright- \ve(x,y) \bleft y, (x\circ z)\circ t\bright  \\
   &&+\ve(x+y,z)\bleft z, (x\circ  y)\circ t\bright-0- (\bleft x, y\bright\circ z)\circ t\\
   &&+ \ve(y,z)(\bleft x, z\bright\circ y)\circ t-\ve(x,y+z)(y\circ z)\circ\bleft x,t\bright\\
   && - (x\circ\bleft y, z\bright)\circ t + \ve(y,z)(x\circ z)\bleft y,t\bright- (x\circ y)\circ\bleft z,t\bright\Big)\\
   &=&-\eight\Big(x\circ [(y\circ z)\circ t]- \ve(x,y)  y\circ [(x\circ z)\circ t]  \\
   &&+\ve(x+y,z) z\circ[(x\circ  y)\circ t]-[(x\circ y)\circ z]\circ t+\ve(x,y)[(y\circ x)\circ z]\circ t\\
   &&+ \ve(y,z)[(x\circ z)\circ y]\circ t-\ve(y,z)\ve(x,z)[(z\circ x)\circ y]\circ t  -\ve(x,y+z)(y\circ z)\circ(x\circ t)\\
   &&+ [x\circ(y\circ z)]\circ t - \ve(y,z)[x\circ(z\circ y)]\circ t + \ve(y,z)(x\circ z)\circ (y\circ t)- (x\circ y)\circ (z\circ t)\Big)\\
   &=&-\eight\Big([x\circ (y\circ z)]\circ t+ \underline{\ve(x,y+z) (y\circ z)\circ (x\circ t)}\\
   &&- \ve(x,y)  [y\circ (x\circ z)]\circ t \underbrace{- \ve(x,y)\ve(y,x+z)(x\circ z)\circ (y\circ t)}  \\
   &&+\ve(x+y,z) [z\circ(x\circ  y)]\circ t+\underline{\underline{\ve(x+y,z)\ve(z,x+y) (x\circ  y)\circ (z\circ t)}}\\
   &&-[(x\circ y)\circ z]\circ t+\ve(x,y)[(y\circ x)\circ z]\circ t\\
   &&+ \ve(y,z)[(x\circ z)\circ y]\circ t-\ve(y,z)\ve(x,z)[(z\circ x)\circ y]\circ t  \underline{-\ve(x,y+z)(y\circ z)\circ(x\circ t)}\\
   &&+ [x\circ(y\circ z)]\circ t - \ve(y,z)[x\circ(z\circ y)]\circ t +
   \underbrace{\ve(y,z)(x\circ z)\circ (y\circ t)}\underline{\underline{-(x\circ y)\circ (z\circ t)}}\Big)\\
   &=&-\textstyle{\frac{3}{8}}\Big(x\circ (y\circ z)-(x\circ y)\circ z-\ve(x,y)y\circ (x\circ z)\Big)\circ t\\
   &=&0,
\end{eqnarray*}
where the underline terms are canceled out.
The other cases can be checked similarly.
The proof is completed.
\qed

\begin{Example}\label{example3}
Recall that in Example \ref{example1}  the $\mathbb{Z}_2 \times \mathbb{Z}_2$-graded analogue of  $\mathfrak{sl}_2(\C)$ consists of  $\g_{(0,1)} = \mathbb{C}x$, $\g_{(1,0)} = \mathbb{C}y$ and $\g_{(1,1)} = \mathbb{C}z$
%%%%%%%%%%%%%%%%%%%%%%%%%%%%%%
\emptycomment{
where the basis elements $x,y,z$ are given by
$$x=\left(\begin{array}{ccc}
0 & 0 & 0 \\ 0 & 0 & -1\\ 0 & 1 & 0\\
\end{array} \right),\quad
y=\left(\begin{array}{ccc}
0 & 0& 1 \\0 & 0 & 0\\-1 &  0 & 0\\
\end{array} \right),\quad
z=\left(\begin{array}{ccc}
0 & -1 & 0 \\1 & 0 & 0\\0 &  0 & 0\\
\end{array} \right).$$
}
%%%%%%%%%%%%%%%%%%%%%
where the bracket is given by
\begin{eqnarray*}
&&[x,x]=[y,y]=[z,z]=0,\\
&&[x, y]=xy+y x=z,\\
&& [y, z]=yz+z y=x,\\
&& [z, x]=zx+xz=y.
\end{eqnarray*}
Let $V=V_{(0,1)}\oplus V_{(1,0)}\oplus V_{(1,1)}=\mathbb{C}e_1\oplus \mathbb{C}e_2\oplus \mathbb{C}e_3$ be the three-dimensional representation
%%%%%%%%%%%%%%%%%%%%%%%%%%%%%%
\emptycomment{
$$e_1=\left(
       \begin{array}{c}
         1 \\
         0 \\
         0 \\
       \end{array}
     \right)
, \quad e_2=\left(
       \begin{array}{c}
         0 \\
         1 \\
         0 \\
       \end{array}
     \right), \quad e_3=\left(
       \begin{array}{c}
         0 \\
         0 \\
         1 \\
       \end{array}
     \right)$$
     }
%%%%%%%%%%%%%%%%%%%
where the action of $\mathfrak{sl}_2^c$ is given by:
\begin{eqnarray*}
xe_1  = 0, \quad ye_1 = -e_3, \quad  ze_1  = e_2 , \\
xe_2  = e_3, \quad ye_2 = 0, \quad  ze_2  = -e_1 , \\
xe_3  = -e_2, \quad ye_3 = e_1, \quad  ze_3  = 0.
\end{eqnarray*}
Then we get a Leibniz color algebra $\mathfrak{sl}_2^c\ltimes_D V$ where $V$ is equal to its left center.
Thus we obtain a 2-term color $L_{\infty}$-algebra as follows:
$$(V\stackrel{i}{\hookrightarrow}\mathfrak{sl}_2^c\ltimes_D V,\, l_2, \,l_3)$$
where $l_2$ is equal to the bracket defined as above on $\mathfrak{sl}_2^c$, $l_2$ is zero on $V$, the other case is given by
\begin{eqnarray*}
l_2(x,e_1) = 0, \quad l_2(y,e_1) = -\half e_3, \quad  l_2(z,e_1)  = \half e_2 , \\
l_2(x,e_2) = \half e_3, \quad l_2(y,e_2) = 0, \quad l_2( z,e_2)  = -\half e_1 , \\
l_2(x,e_3)  = -\half e_2, \quad l_2(y,e_3 )= \half e_1, \quad  l_2(z,e_3)  = 0.
\end{eqnarray*}
By direct computations, $l_3$ is given by
\begin{eqnarray*}
&&l_3(x,  y, e_1)=-\four e_2,  \quad l_3(x,  y, e_2)=\four e_1, \quad l_3(x,  y, e_3)=0,\\
&&l_3(y,  z, e_1)=0,  \quad l_3(y,  z, e_2)=-\four e_3, \quad l_3(y,  z, e_3)=\four e_2,\\
&&l_3(x,  z, e_1)=\four e_3,  \quad l_3(x,  z, e_2)=0, \quad l_3(x,  z, e_3)=-\four e_1.
\end{eqnarray*}
\end{Example}

\section{Acknowledgements}

The research was supported by NSFC (11501179) and a doctoral research program (qd14148) of Henan Normal University.
I would like to thank Professor Dimitry Leites for sending me the book \cite{Lei11}.
I also thank Professor Zhangju Liu and Professor Yunhe Sheng for their helpful suggestions and discussions.
%Special thanks to the anonymous  referee for valuable comments which helped to improve this paper significantly.

\vskip7pt

\footnotesize{
\noindent College of Mathematics and Information Science,\\
Henan Normal University, Xinxiang 453007, P. R. China;\\
 E-mail address:\texttt{{  zhangtao@htu.cn}}
}

\end{document}